\documentclass [12pt,leqno] {article}
\usepackage{amsmath,amsfonts,url}
\usepackage{enumerate}

\def\supp{\mathop{\rm supp}\nolimits}

\newcommand{\dowod}{{\em Proof}.\/ }
\newcommand{\qed}{~\hfill~$\fbox{}$ \vspace{0.5cm}}
\newcommand{\R}{ \mathbb{R}^{d}}
\newcommand{\N}{\mathbb{N}}
\newcommand{\indyk}[1]{{\bf 1}_{#1}}
\newcommand{\sfera}{ \mathbb{S}}
\newcommand{\gener}{{\cal A}}
\newcommand{\scalp}[2]{#1\cdot#2}

\newcommand{\supmasyL}{\bar{b}_\varepsilon}
\newcommand{\infmasyL}{\underline{b}_\varepsilon}
\newcommand{\Cvanish}{C_\infty(\R)}
\newcommand{\Bbound}{B_b(\R)}
\newcommand{\Aapprox}[1]{{\cal A}_{#1}}

\newtheorem{lemma}{\indent\sc Lemma}
\newtheorem{theorem}{\indent\sc Theorem}

\begin{document}

\title{Approximation of Stable-dominated Semigroups}
\author{Pawe{\l} Sztonyk}
\footnotetext{ Institute of Mathematics and Computer Science,
  Wroc{\l}aw University of Technology,
  Wybrze{\.z}e Wyspia{\'n}\-skie\-go 27,
  50-370 Wroc{\l}aw, Poland.\\
  {\rm e-mail: Pawel.Sztonyk@pwr.wroc.pl} \\
}
\maketitle

\begin{abstract}
  We consider Feller semigroups with jump intensity dominated by that of the rotation invariant stable L\'evy    process.
  Using an approximation scheme  
  we obtain estimates of corresponding heat kernels.  
\end{abstract}
\footnotetext{2000 {\it MS Classification}:
Primary 60J75, 60J35; Secondary 47D03.
}

\section{Introduction and Preliminaries}\label{Prel}

Let $p(t,x,y)$ be the transition density  of the rotation invariant $\alpha$--stable L\'evy process on $\R$ with the L\'evy measure
\begin{equation}\label{eq:RISLM}
  \nu(dy)=\frac{dy}{|y|^{\alpha+d}},\quad y\in\R\setminus\{0\}.
\end{equation}
Here $\alpha\in(0,2)$ and $d=1,2,\dots$
It is well-known (see, e.g., \cite{BlGe}) that
$$
  p(t,x,y) \approx \min\left(t^{-d/\alpha},\frac{t}{|y-x|^{\alpha+d}}\right),\quad t>0,\,x,y\in\R.
$$

The upper bound $p(t,x,y)\leq c \min\{t^{-d/\alpha},t|y-x|^{-\alpha-d}\}$, holds for
every symmetric stable L\'evy process whose L\'evy measure is bounded above by a constant multiple of $\nu$ (see \cite{GH}). 
More diverse asymptotic results for stable L\'evy processes were given 
in papers \cite{W} and \cite{BS2007}. In particular if for some $\gamma\in[1,d]$ 
the L\'evy measure $\bar{\nu}$ of the stable process satisfies
$$
  \bar{\nu}(B(x,r))\leq c r^{\gamma} \quad  \mbox{for all}\quad |x|=1,\, r\leq 1/2,
$$  
then we have 
$$
  \bar{p}(t,x,y)\leq c \min\left(t^{-d/\alpha},t^{1+\frac{\gamma-d}{\alpha}}
                     |y-x|^{-\alpha-\gamma}\right),\quad x,y\in\R,
$$
where $\bar{p}$ denotes the corresponding transition density.

Similar results hold for other L\'evy and Markov processes. Estimates of transition
densities for tempered stable processes are given in \cite{Szt2008}. The case of stable--like and
mixed type Markov processes was investigated in \cite{ChKum} and \cite{ChKum08} by Z.-Q.~Chen and T.~Kumagai. 
K.~Bogdan and
T.~Jakubowski in \cite{BJ2007} obtained estimates of heat kernels of the fractional Laplacian perturbed by gradient operators. Derivatives of 
stable densities have been considered in \cite{Lewand} and \cite{Szt2007}.

Our main goal are the heat kernel estimates for a class of stable--dominated Feller semigroups.
The name {\it stable-- dominated} refers to the inequality (A.1).

Let $f:\R\times\R\mapsto [0,\infty]$ be a Borel function. We consider the following conditions on $f$. \\

\textbf{(A.1)} There exists a constant $M>0$ such that
$$
  f(x,y) \leq \frac{M}{|y-x|^{\alpha+d}},\quad x,y\in\R,\, y\neq x.
$$
\\

\textbf{(A.2)}  $\,f(x,x+h)=f(x,x-h)\,$ for all $x,h\in\R$, or $\alpha<1$.\\

\textbf{(A.3)} $f(x,y)=f(y,x)$ for all $x,y\in\R$.\\

\textbf{(A.4)} There exists $a>0$ such that
$$
  \inf_{x\in\R} \int_{|y-x|>\varepsilon} f(x,y)\,dy \geq a \varepsilon^{-\alpha},\quad \varepsilon>0.
$$

\noindent We denote 
$$
  b_\varepsilon(x)=\int_{|y-x|>\varepsilon} f(x,y)\,dy,\,\quad \varepsilon>0,\,x\in\R.
$$
It follows from (A.1) that there exists a constant $A>0$ such that
\begin{equation}\label{eq:supbx}
  \supmasyL := \sup_{x\in\R}b_\varepsilon(x) \leq A\varepsilon^{-\alpha},\quad \varepsilon>0.
\end{equation}
Thus, (A.4) is a partial converse of (A.1) and we have
$$
  \infmasyL := \inf_{x\in\R}b_\varepsilon(x) \geq a \varepsilon^{-\alpha},\quad \varepsilon>0.
$$

For $x\in\R$ and $r>0$ we let $B(x,r)=\{y\in\R :\:|y-x|<r\}$. $B_b(\R)$ denotes the set of bounded Borel measurable functions, $C^k_c(\R)$ denotes the set of $k$ times continuously differentiable functions with compact support and $\Cvanish$ is the set of 
continuous functions vanishing at infinity. 
We use $c,C$ (with subscripts) to denote finite positive constants
which depend only on the constants $\alpha,M,a$, and the dimension $d$. Any {\it additional} dependence
is explicitly indicated by writing, e.g., $c=c(n)$.
The value of $c,C$, when used without subscripts, may
change from place to place. We write $f(x)\approx g(x)$ to indicate that there is a constant $c$ such that $c^{-1}f(x) \leq g(x) \leq c f(x)$. 

Assuming (A.1) and (A.2) we may consider the operator
\begin{eqnarray*}
\gener\varphi(x) & =  &
\lim_{\varepsilon\downarrow 0}\int_{|y-x|>\varepsilon} \left(\varphi(y)-\varphi(x)\right)f(x,y)\,dy\\
&  =   & \int_{\R} \left(\varphi(x+h)-\varphi(x)-\scalp{h}{\nabla\varphi(x)}\indyk{|h|<1}\right)
          f(x,x+h)\,dh \\
&      & +\, \frac{1}{2}\int_{|h|<1} \scalp{h}{\nabla\varphi(x)} \left(f(x,x+h)-f(x,x-h) \right)\,dh ,\quad \varphi\in C^2_c(\R).
\end{eqnarray*}

We record the following basic fact (we postpone the proof to Section \ref{Est}).
\begin{lemma}\label{l:AonCvanish}
If (A.1), (A.2) hold and the function
$x\to f(x,y)$ is continuous on $\R\setminus\{y\}$ for every $y\in\R$ then $\gener$ maps $C^2_c(\R)$ into $\Cvanish$.
\end{lemma}

In the following we always assume that the condition (A.1) is satisfied.
For every $\varepsilon>0$ we denote
$$
  f_{\varepsilon}(x,y)=\indyk{B(0,\varepsilon)^c}(y-x)f(x,y),\quad x,y\in\R,
$$
and
$$
  \Aapprox{\varepsilon} \varphi(x)=\int \left(\varphi(y)-\varphi(x)\right)f_{\varepsilon}(x,y)\,dy,\quad \varphi\in \Bbound.
$$

We note that the operators $\Aapprox{\varepsilon}$ are bounded since 
$|\Aapprox{\varepsilon}\varphi(x)| \leq 2 \| \varphi \|_\infty b_\varepsilon(x)\leq 2\| \varphi \|_\infty \supmasyL$.
Therefore the operator 
$$
  e^{t\Aapprox{\varepsilon}}=\sum_{n=0}^\infty \frac{t^n\Aapprox{\varepsilon}^n}{n!},\quad t\geq 0,\,\varepsilon>0,
$$
is well--defined and bounded: $\Bbound\to\Bbound$. In fact for every $\varepsilon>0$ the family of operators $\{e^{t\Aapprox{\varepsilon}},\,t\leq 0\}$ 
is a semigroup on $\Bbound$, i.e., $e^{(t+s)\Aapprox{\varepsilon}}=e^{t\Aapprox{\varepsilon}}e^{s\Aapprox{\varepsilon}}$ 
for all $t,s\geq 0$, $\varphi\in\Bbound$. We note that $e^{t\Aapprox{\varepsilon}}$ is positive for all $t\geq 0$, $\varepsilon>0$ (see (\ref{eq:exp})).

Our main result is the following theorem.

\begin{theorem}\label{th:MainT} 
If (A.1) -- (A.4) are satisfied then there exists a constant $C$ such that 
for every nonnegative $\varphi\in \Bbound$ and $\varepsilon>0$ we have
\begin{equation}\label{eq:MainT}
  e^{t\Aapprox{\varepsilon}}\varphi(x) \leq C \int 
  \varphi(y) \min\left(t^{-d/\alpha},\frac{t}{|y-x|^{\alpha+d}}\right) dy + e^{-tb_\varepsilon(x)}\varphi(x),\quad x\in\R.
\end{equation}
\end{theorem}

We will prove Theorem \ref{th:MainT} in Section \ref{Est} after a sequence of lemmas. To study a limiting semigroup 
we will use additional assumptions.\\

\textbf{(A.5)} The function
$x\to f(x,y)$ is continuous on $\R\setminus\{y\}$ for every $y\in\R$.\\

\textbf{(A.6)} $\gener$ regarded as an operator on $\Cvanish$ is closable and its closure $\bar{\gener}$ is a generator of a strongly continuous 
contraction semigroup of operators $\{P_t,\,t\geq 0\}$ on $\Cvanish$.\\

We note that the operator $\gener$ satisfies {\it the positive maximum principle}, i.e., if for 
$\varphi\in C^2_c(\R)$ we have $\sup_{x\in\R}\varphi(x)=\varphi(x_0)\geq 0$ 
then $\gener\varphi(x_0)\leq 0$. This yields that $P_t$ is positive operator for each $t\geq0$ (see Theorem 1.2.12 and Theorem 4.2.2 in \cite{EK}). 
Thus, by our assumptions, $\{P_t,\,t\geq 0\}$ is {\it a Feller semigroup}.

In Section \ref{Est} we prove the following theorem.

\begin{theorem}\label{th:main}
If (A.1)--(A.6) hold then there is $p:(0,\infty)\times\R\times\R\to [0,\infty)$
 such that
$$
  P_t \varphi (x) = \int_{\R} \varphi(y) p(t,x,y)\, dy,\quad x\in\R,\,t>0,\,\varphi\in \Cvanish,
$$
and 
\begin{equation}\label{eq:main}
  p(t,x,y) \leq C \min\left(t^{-d/\alpha}, \frac{t}{|y-x|^{\alpha+d}}\right),\quad x,y\in\R,\,t>0,
\end{equation}
for some constant $C>0$.
\end{theorem}

We note that $\gener$ is conservative meaning that 
for $\phi\in C^{\infty}_c(\R)$ such that $0\leq \phi\leq 1$, $\phi(0)=1$, and 
$\phi_k(x)=\phi(x/k)$, we have
$\sup_{k\in\N}\|\gener\phi_k\|_\infty<\infty$, and 
$\lim_{k\to\infty}(\gener\phi_k)(x)=0$, for every
$x\in\R$. It follows from Theorem 4.2.7 in \cite{EK} that 
there exists a Markov process $\{X_t,\,t\geq 0\}$ such that
${\mathbb E}[\varphi(X_t)|X_0=x]=P_t\varphi(x)$.

We recall that every generator $G$ of a Feller semigroup, such that $C_c^\infty(\R)\subset {\cal D}(G)$
is necessarily of the form
\begin{eqnarray}\label{eq:FellerGen}
  G\varphi(x)
  &  =  & \sum_{i,j=1}^d q_{ij}(x)D_{x_i}D_{x_j} \varphi(x) + l(x) \nabla \varphi(x) -c(x)\varphi(x) \\
  &     & + \int\limits_{\R}\left(\varphi(x+h)-\varphi(x)-
            \scalp{h}{\nabla \varphi(x)}\;\indyk{|h|<1}\right)\, \nu(x,dh)\,. \nonumber
\end{eqnarray}
Here $\varphi\in C_c^\infty(\R)$, $q(x)=(q_{ij}(x))_{i,j=1}^n$ is a nonnegative definite real symmetric matrix,
the vector $l(x)=(l_i(x))_{i=1}^d$ has real coordinates, $c(x)\geq 0$,
and $\nu(x,\cdot)$ is a L\'evy measure. 
This description is due to Courr\'ege, see \cite[Chapter 4.5]{Jc1}.

The problem whether a given operator $G$ generates a Feller semigroup is not completely resolved yet.
For the interested reader we remark that criteria are given, 
e.g., in \cite{Hoh1,Hoh2,Jacob1,Jacob2,Jacob3}. Generally, smoothness of the coefficients $q,l,c,\nu$
in (\ref{eq:FellerGen}) is sufficient for the existence .
In particular the next lemma follows from Theorem 5.24 in \cite{HohHab} (see also Theorem 4.6.7 in \cite{Jacob4}). 
We omit the proof as it is a straightforward
verification of the assumptions given there. 
\begin{lemma}Let $m$ be a smallest integer such that $m>\max(\frac{d}{\alpha},2)+d$ and define
$k=2m+1-d$. Assume that
  \begin{enumerate}[(i)]
  	\item there exist a cone $\emptyset \neq V\subset\R$ and a constant $c$ such that 
      $$
        c^{-1}\indyk{V}(h) |h|^{-\alpha-d}\leq f(x,x+h)\leq c|h|^{-\alpha-d},\quad x\in\R,\, h\in\R\setminus\{0\},
      $$
    \item 
      $$
        f(x,x+h)=f(x,x-h),\quad x,h\in\R,
      $$
    \item for every $h\in\R\setminus\{0\}$ the function $g_h(x)=f(x,x+h)$ is 
      $k$-times continuously differentiable and
      $$
        |\partial_x^\beta g_h(x)| \leq c |h|^{-\alpha-d},\quad x\in\R,\,h\in\R\setminus\{0\},
      $$
      for each $\beta\in\N_0^d,$ $|\beta|\leq k$.
  \end{enumerate}
  Then $\gener$ has an extension that generates a Feller semigroup, i.e., the condition (A.6) is satisfied.
\end{lemma}

For the stable -- like case more recent results are given
by R.~Schilling and T.~Uemura in \cite{SchUe}. In particular it follows from Corollary 6.4 in \cite{SchUe} that
if $f(x,y)\approx |y-x|^{-\alpha-d}$ for $|y-x|\leq 1$, 
$0\leq f(x,y)\leq c |y-x|^{-\gamma-d}$ for $|y-x|> 1$ and some $\gamma>0$, and the function $(x,y)\to\log f(x,y)/\log |y-x|$ is Lipschitz continuous,
then $\gener$ generates a Feller semigroup.

Z.-Q.~Chen and T.~Kumagai in \cite{ChKum} and \cite{ChKum08} investigate the case of symmetric jump--type Markov processes 
on metric measure spaces by using the Dirichlet form. They prove
the existence and obtain estimates
of the densities (see Theorem~1.2 in \cite{ChKum08}) analogous to (\ref{eq:main}).
The jump
kernels in \cite{ChKum} and \cite{ChKum08} are assumed to be {\it comparable} with certain rotation invariant functions.
In Theorem \ref{th:main} we assume the estimate (A.1) from above but we use (A.4) as the only
estimate for the size of $f$ from below. This is the main novelty as far as the class of jump kernels is considered.
We also propose a new technique of estimating the semigroups $\{ e^{t\Aapprox{\varepsilon}},\,t\geq 0\}$, which may be of more general interest. We note that an alternative approximation scheme
is given in \cite{BSchill}. We also note that the estimate (\ref{eq:main}) can not be generally improved as seen in the case 
of L\'evy processes 
(see \cite{BS2007}, \cite{W}, \cite{Szt2008}). 
In fact, if $\bar{\nu}(dy)=g(y/|y|)|y|^{-\alpha-d}\,dy$ and $g$ is continuous on the unit sphere $\sfera$ 
we have
$\lim_{r\to\infty}r^{d+\alpha}p(1,0,r\theta)=cg(\theta)$, $\theta\in\sfera$. If $g(\theta)=0$ then $\lim_{r\to\infty}r^{d+2\alpha}p(1,0,r\theta)=c_{\theta}>0$, as proved by J. Dziuba\'nski in \cite{Dziub}.


\section{Approximation}\label{Est}

We first observe that
\begin{eqnarray}\label{eq:KGamma}
  \Aapprox{\varepsilon} \varphi (x)  
  &  =  &  \int \left(\varphi(y)-\varphi(x)\right)f_{\varepsilon}(x,y)\,dy + 
           (\supmasyL-b_\varepsilon(x))\int (\varphi(y)-\varphi(x))\delta_x(dy) \nonumber \\
  &  =  &  \int (\varphi(y)-\varphi(x))\tilde{\nu}_{\varepsilon}(x,dy)  \\
  &  =  & \Gamma_\varepsilon \varphi(x) - \supmasyL\varphi(x),\quad \varphi\in \Cvanish,\,x\in\R,\nonumber        
\end{eqnarray}
where
$$
  \tilde{\nu}_\varepsilon(x,dy)=f_{\varepsilon}(x,y)\,dy+(\supmasyL-b_\varepsilon(x))\delta_x(dy),
$$
and
$$
  \Gamma_\varepsilon \varphi(x)=\int \varphi(y)\tilde{\nu}_\varepsilon(x,dy),\quad \varphi\in \Cvanish,x\in\R.
$$
This yields that
\begin{equation}\label{eq:exp}
  e^{t\Aapprox{\varepsilon}}\varphi(x)=e^{t(\Gamma_\varepsilon-\supmasyL I)}\varphi(x)=
  e^{-t\supmasyL}e^{t\Gamma_\varepsilon}\varphi(x).
\end{equation}
A consequence of (\ref{eq:exp}) is that we may consider the operator $\Gamma_\varepsilon$ 
and its powers instead of $\Aapprox{\varepsilon}$.
The fact that $\Gamma_\varepsilon$ is positive enables for more precise estimates.

For $n\in\N$ we define
\begin{eqnarray*}
  f_{n+1,\varepsilon}(x,y)
  &  =   &  \int f_{n,\varepsilon}(x,z) f_{\varepsilon}(z,y)\,dz  \\
  &      &  +\, \left(\supmasyL-b_\varepsilon(y)\right)f_{n,\varepsilon}(x,y) 
            \,+\, \left(\supmasyL-b_\varepsilon(x)\right)^{n}f_{\varepsilon}(x,y),
\end{eqnarray*}
where we let $f_{1,\varepsilon}=f_\varepsilon$.
By induction and Fubini--Tonelli theorem we get
\begin{equation}\label{eq:intf_n}
    \int f_{n,\varepsilon}(x,y)\,dy = \supmasyL^n-\left(\supmasyL-b_\varepsilon(x)\right)^n\,,\quad x\in\R,\,n\in\N.
\end{equation}

\begin{lemma}\label{lm:GammaDensity} For all $\varepsilon>0$, $x\in\R$, and $n\in\N$ we have
\begin{equation}\label{eq:GammaDensity}
  \Gamma_\varepsilon^n \varphi(x) = \int \varphi(z) f_{n,\varepsilon}(x,z)\,dz 
  + \left(\supmasyL-b_\varepsilon(x)\right)^{n}\varphi(x),\quad \varphi\in \Cvanish.
\end{equation}
\end{lemma}
\dowod We use induction. For $n=1$ we have (\ref{eq:GammaDensity}) from the definition of $\Gamma_\varepsilon$. Let us assume that (\ref{eq:GammaDensity}) holds for some $n\in\N$.
Using Fubini's theorem we get
\begin{eqnarray*}
  \Gamma^{n+1}_\varepsilon\varphi(x) 
  &  =  & \int \Gamma_\varepsilon \varphi(z) f_{n,\varepsilon}(x,z)\,dz +
          \left(\supmasyL-b_\varepsilon(x)\right)^n \Gamma_\varepsilon \varphi(x) \\
  &  =  & \int \left(\int \varphi(y)f_{\varepsilon}(z,y)\,dy +
          \left(\supmasyL-b_\varepsilon(z)\right) \varphi(z)\right)f_{n,\varepsilon}(x,z)\,dz \\
  &     & + \left(\supmasyL-b_\varepsilon(x)\right)^n \left(\int \varphi(y)f_{\varepsilon}(x,y)\,dy +
          \left(\supmasyL-b_\varepsilon(x)\right)\varphi(x)\right) \\
  &  =  & \int \varphi(z) f_{n+1,\varepsilon}(x,z)\,dz + 
          \left(\supmasyL-b_\varepsilon(x)\right)^{n+1}\varphi(x).
\end{eqnarray*}
\qed

The following lemma is a key observation in our development. The significance 
of the inequality (\ref{eq:LI}) is that on the right hand side we obtain precisely $b_\varepsilon(y) |y-x|^{-\alpha-d}$.
One can interpret (\ref{eq:LI}) as subharmonicity of
the function $\varphi(z)=\indyk{B(y,\kappa|y-x|)}(z)|z-x|^{-\alpha-d}$ at $y$ with respect to all operators
$\Aapprox{\varepsilon}$.

\begin{lemma}\label{lm:EintL}
  If (A.1), (A.2) and (A.4) hold then there exists $\kappa\in(0,1)$ such that
  \begin{equation}\label{eq:LI}
    \int_{B(y,\kappa|y-x|)} |z-x|^{-\alpha-d}f_{\varepsilon}(y,z)\,dz
    \leq b_\varepsilon(y) |y-x|^{-\alpha-d}
  \end{equation}
  for every $x,y\in\R,\, \varepsilon>0$.
\end{lemma}
\dowod We fix $x,y\in\R$ and let $\phi(z)=|y-x|^{\alpha+d}|z-x|^{-\alpha-d}$, $|z-x|>0$. We have $\phi(y)=1$, 
\begin{equation}\label{eq:1derphi}
  \partial_j \phi(z)=(\alpha+d)|y-x|^{\alpha+d}|z-x|^{-\alpha-d-2}(x_j-z_j),
\end{equation}
and
\begin{displaymath}
  \partial_{j,k} \phi(z) = (\alpha+d)|y-x|^{\alpha+d}|z-x|^{-\alpha-d-2}\left[(\alpha+d+2)\frac{(x_j-z_j)(x_k-z_k)}{|x-z|^{2}}-\delta_{jk}\right].
\end{displaymath}
This yields
\begin{equation}\label{eq:2derphi}
  \sup_{\substack{ z\in B(y,\kappa|y-x|),\\ j,k\in\{1,\dots,d\} }}  |\partial_{j,k} \phi(z)| \leq (\alpha+d)(\alpha+d+3)(1-\kappa)^{-\alpha-d-2}|y-x|^{-2},
\end{equation}
for every $\kappa\in(0,1)$.
Using the Taylor expansion for $\phi$, (\ref{eq:1derphi}) and (\ref{eq:2derphi}) in the second inequality below,
(A.1) and (A.2) in the third, and (A.4) in the forth we obtain
\begin{eqnarray*}
  &      & \left| \int_{B(y,\kappa|y-x|)}\left[ 1 - |y-x|^{\alpha+d}|z-x|^{-\alpha-d}\right]f_{\varepsilon}(y,z)\,dz\right| \\
  &   =  & \left| \int_{B(0,\kappa|y-x|)} 
                  \left(\phi(y)-\phi(y+h)\right) f_{\varepsilon}(y,y+h)\,dh \right| \\
  & \leq & \left| \int_{B(0,\kappa|y-x|)} 
                  \left(\scalp{\nabla\phi (y)}{h}+ \phi(y)-\phi(y+h)\right)
                   f_{\varepsilon}(y,y+h)\,dh \right| \\
  &      & + \left| \int_{B(0,\kappa|y-x|)} \scalp{\nabla\phi (y)}{h}
                   \frac{f_{\varepsilon}(y,y-h)-f_{\varepsilon}(y,y+h)}{2}\,dh \right| \\
  & \leq & c_1(1-\kappa)^{-\alpha-d-2} |y-x|^{-2} 
           \int_{B(0,\kappa|y-x|)} |h|^2 f(y,y+h)\,dh \\
  &      & + c_2 |y-x|^{-1} \int_{B(0,\kappa|y-x|)} |h|
                   \left|f_{\varepsilon}(y,y-h)-f_{\varepsilon}(y,y+h)\right|\,dh  \\
  & \leq & c_3|y-x|^{-\alpha}\kappa^{1-\alpha} \left(\kappa(1-\kappa)^{-\alpha-d-2}+1\right)
   \leq  b_{\kappa|y-x|}(y),
\end{eqnarray*}
for sufficiently small $\kappa=\kappa(\alpha,d)$.
This yields
\begin{eqnarray*}
   &      & \int_{B(y,\kappa|y-x|)} |z-x|^{-\alpha-d}f_{\varepsilon}(y,z)\,dz \\
   & \leq & \left| \int_{B(y,\kappa|y-x|)}
          \left[ 1 - |y-x|^{\alpha+d}|z-x|^{-\alpha-d}\right]f_{\varepsilon}(y,z)\,dz\right||y-x|^{-\alpha-d} \\
   &      & + \int_{B(y,\kappa|y-x|)} f_{\varepsilon}(y,z)\,dz\, |y-x|^{-\alpha-d} \\
   & \leq & b_\varepsilon(y) |y-x|^{-\alpha-d}. 
\end{eqnarray*}
\qed

We can now obtain estimates of $f_{n,\varepsilon}(x,y)$. The first one is an adequate description of the decay rate
at infinity while the next two give global bounds.

\begin{lemma}\label{lm:Estimate1}
  If (A.1) -- (A.4) hold then there exists a constant $C$ such that
  \begin{equation}\label{eq:Estimate1}
    f_{n,\varepsilon}(x,y) \leq C n\supmasyL^{n-1}|y-x|^{-\alpha-d},\quad 
    x,y\in\R,\,\varepsilon>0,\,n\in\N.
  \end{equation}
\end{lemma}
\dowod We use induction. For $n=1$ the inequality (\ref{eq:Estimate1}) holds with $C=M$. Let $\kappa\in (0,1)$ be such that (\ref{eq:LI}) is satisfied. 
We will prove that (\ref{eq:Estimate1}) holds
with $C=M\kappa^{-\alpha-d}$. We have
$$
  \int f_{n,\varepsilon}(x,z)f_{\varepsilon}(z,y)\, dz = \int_{B(y,\kappa|y-x|)^c} + \int_{B(y,\kappa|y-x|)}
  = I+II.
$$
By (A.1) and (\ref{eq:intf_n}) we get
\begin{eqnarray*}
  I
  &   =  & \int_{B(y,\kappa|y-x|)^c} f_{n,\varepsilon}(x,z) f_{\varepsilon}(z,y)\, dz \\
  & \leq & \int_{B(y,\kappa|y-x|)^c} f_{n,\varepsilon}(x,z)\frac{M}{|y-z|^{\alpha+d}}\, dz \\
  & \leq & \frac{M}{\kappa^{\alpha+d}|y-x|^{\alpha+d}} \int f_{n,\varepsilon}(x,z)\,dz \\
  &   =  & \frac{M}{\kappa^{\alpha+d}|y-x|^{\alpha+d}} 
           \left[\supmasyL^n-\left(\supmasyL-b_\varepsilon(x)\right)^n\right].
\end{eqnarray*}
By the symmetry of $f$ (see (A.3)), induction and Lemma \ref{lm:EintL} we obtain
\begin{eqnarray*}
  II
  &  =   &  \int_{B(y,\kappa|y-x|)} f_{n,\varepsilon}(x,z) f_{\varepsilon}(z,y)\, dz \\
  & \leq & C n \supmasyL^{n-1} \int_{B(y,\kappa|y-x|)}
           |z-x|^{-\alpha-d} f_{\varepsilon}(z,y)\,dz \\
  &  =   & C n \supmasyL^{n-1} \int_{B(y,\kappa|y-x|)}
           |z-x|^{-\alpha-d} f_{\varepsilon}(y,z)\,dz \\
  & \leq & C n \supmasyL^{n-1} b_\varepsilon(y) |y-x|^{-\alpha-d} .
\end{eqnarray*}
We get
\begin{eqnarray*}
  f_{n+1,\varepsilon}(x,y)
  &  =   & I+II + \left(\supmasyL-b_\varepsilon(y)\right)f_{n,\varepsilon}(x,y) +
          \left(\supmasyL-b_\varepsilon(x)\right)^{n}f_{\varepsilon}(x,y) \\
  & \leq & M \kappa^{-\alpha-d}|y-x|^{-\alpha-d} 
           \left[\supmasyL^n-\left(\supmasyL-b_\varepsilon(x)\right)^n\right]
           + C n \supmasyL^{n-1} b_\varepsilon(y) |y-x|^{-\alpha-d}  \\
  &      & + \left(\supmasyL-b_\varepsilon(y)\right)C n \supmasyL^{n-1} |y-x|^{-\alpha-d} +
          \left(\supmasyL-b_\varepsilon(x)\right)^{n}M|y-x|^{-\alpha-d} \\
  &   <  & C (n+1) \supmasyL^{n} |y-x|^{-\alpha-d}.
\end{eqnarray*}
\qed

\begin{lemma}\label{lm:Estimate2} Assume (A.1), (A.3) and (A.4). Then
there exists $C$ such that
\begin{equation}\label{eq:Estimate2}
    f_{n,\varepsilon}(x,y)\leq C \supmasyL^{d/\alpha}\left(\supmasyL^n-\left(\supmasyL-b_\varepsilon(x)\right)^n\right),
    \quad x,y\in\R,\,\varepsilon>0,\,n\in\N.
  \end{equation}  
\end{lemma}
\dowod For $n=1$ by (A.1) and (A.4) we have 
$$
  f_{\varepsilon}(x,y)\leq M \varepsilon^{-\alpha-d} \leq 
  M\left( \frac{b_\varepsilon(x)}{a} \right)^{(\alpha+d)/\alpha}
  \leq
  M 
  \left(\frac{b_\varepsilon(x)}{a} \right)
  \left(\frac{\supmasyL}{a}\right)^{d/\alpha}  ,
$$ 
and so (\ref{eq:Estimate2}) holds with $C=Ma^{-d/\alpha-1}$.
Let (\ref{eq:Estimate2}) holds for some $n\in\N$ with $C=Ma^{-d/\alpha-1}$. 
By induction and the symmetry of 
$f_{\varepsilon}$ we get
\begin{eqnarray*}
  f_{n+1,\varepsilon}(x,y)
  & \leq & C \supmasyL^{d/\alpha}\left(\supmasyL^n-\left(\supmasyL-b_\varepsilon(x)\right)^n\right)\left(\int f_{\varepsilon}(y,z)\,dz + \supmasyL-b_\varepsilon(y)\right)\\
  &      & + \left(\supmasyL-b_\varepsilon(x)\right)^n C \supmasyL^{d/\alpha}b_\varepsilon(x) \\
  &   =  & C (\supmasyL)^{d/\alpha}\left(\supmasyL^{n+1}-\left(\supmasyL-b_\varepsilon(x)\right)^{n+1}\right).
\end{eqnarray*}
\qed

\begin{lemma}\label{lm:Estimate3} If (A.1), (A.3) and (A.4) are satisfied then there exists $C$ such that
\begin{equation}\label{eq:Estimate3}
    f_{n,\varepsilon}(x,y)\leq C \supmasyL^{n+d/\alpha}
    n^{-d/\alpha},
    \quad x,y\in\R,\,\varepsilon>0,\,n\in\N.
  \end{equation}  
\end{lemma}
\dowod 
We may choose $n_0\in\N$ such that 
\begin{equation}\label{eq:defn_0}
  (1-a/A)^n (n+1)^{d/\alpha}<\frac{1}{n+1}
\end{equation} 
for every $n\geq n_0$. For $n\leq n_0$ by Lemma \ref{lm:Estimate2} 
we have 
$$
  f_{n,\varepsilon}(x,y)\leq c_1 \supmasyL^{d/\alpha}\supmasyL^n \leq c_1 \supmasyL^{n+d/\alpha}
    n^{-d/\alpha} n_0^{d/\alpha},
$$
which yields
the inequality (\ref{eq:Estimate3}) with $C=c_1n_0^{d/\alpha}$ in this case. For $n\geq n_0$ we use induction. 
Let 
$$
p=\frac{d\,2^{\max(d/\alpha,1)-1}}{\alpha},\quad \mbox{and}  \quad
  \eta=\left(\frac{(a/A)^2}{2(1+p)}\right)^{\frac{1}{\alpha}}.
$$
We assume
that (\ref{eq:Estimate3}) holds for some $n\geq n_0$ with
$C=\max(c_1n_0^{d/\alpha},M\eta^{-\alpha-d}a^{-1-d/\alpha})$. We have
$$
  \int  f_{n,\varepsilon}(x,z)f_{\varepsilon}(z,y)\,dz = \int_{B(y,\eta \varepsilon(n+1)^{1/\alpha})^c } + \int_{B(y,\eta \varepsilon(n+1)^{1/\alpha})}
  = I+II.
$$
By (A.1), (A.4) and (\ref{eq:intf_n}) we get
\begin{eqnarray*}
  I
  &   =  & \int_{B(y,\eta \varepsilon(n+1)^{1/\alpha})^c} f_{n,\varepsilon}(x,z) f_{\varepsilon}(z,y)\,dz \\
  & \leq & M \int_{B(y,\eta \varepsilon(n+1)^{1/\alpha})^c} f_{n,\varepsilon}(x,z)|y-z|^{-\alpha-d}\,dz \\
  & \leq & M \eta^{-\alpha-d}\varepsilon^{-\alpha-d} (n+1)^{-1-d/\alpha}\int f_{n,\varepsilon}(x,z)\,dz \\
  & \leq & M \eta^{-\alpha-d}a^{-1-d/\alpha}\supmasyL^{1+d/\alpha}(n+1)^{-1-d/\alpha}
           \left[\supmasyL^n-\left(\supmasyL-b_\varepsilon(x)\right)^n\right].
\end{eqnarray*}
By induction, the symmetry of $f_{\varepsilon}$, (\ref{eq:supbx}) and (A.4) we obtain
\begin{eqnarray*}
  II
  &   =  & \int_{B(y,\eta \varepsilon(n+1)^{1/\alpha})} f_{n,\varepsilon}(x,z)f_{\varepsilon}(z,y)\,dz \\
  & \leq & C \supmasyL^{n+d/\alpha}n^{-d/\alpha} \int_{B(y,\eta \varepsilon(n+1)^{1/\alpha})} f_{\varepsilon}(y,z)\,dz \\
  &   =  & C \supmasyL^{n+d/\alpha}n^{-d/\alpha} 
          \left(b_\varepsilon(y)-b_{\eta \varepsilon (n+1)^{1/\alpha}}(y)\right)\\
  & \leq & C \supmasyL^{n+d/\alpha}n^{-d/\alpha} 
          b_\varepsilon(y)\left(1-\frac{a\eta^{-\alpha}}{A(n+1)}\right).
\end{eqnarray*}
By (\ref{eq:defn_0}) we also have 
\begin{equation}\label{eq:hilfe}
  \left(1-\frac{b_\varepsilon(x)}{\supmasyL}\right)^n (n+1)^{d/\alpha} \leq (1-a/A)^n (n+1)^{d/\alpha} \leq
  \frac{1}{n+1}.
\end{equation}
We get
\begin{eqnarray*}
  f_{n+1,\varepsilon}(x,y)
  &   =  & I+II + \left(\supmasyL-b_\varepsilon(y)\right)f_{n,\varepsilon}(x,y) +
          \left(\supmasyL-b_\varepsilon(x)\right)^{n}f_{\varepsilon}(x,y) \\
  & \leq & C \supmasyL^{1+d/\alpha}(n+1)^{-1-d/\alpha}\left[\supmasyL^n-\left(\supmasyL-b_\varepsilon(x)\right)^n\right]\\
  &      & + C \supmasyL^{n+d/\alpha}n^{-d/\alpha} b_\varepsilon(y) 
           \left(1-\frac{a\eta^{-\alpha}}{A(n+1)}\right)\\
  &      & + C \supmasyL^{n+d/\alpha}n^{-d/\alpha} \left(\supmasyL-b_\varepsilon(y)\right)\\
  &      & + C \supmasyL^{1+d/\alpha} \left(\supmasyL-b_\varepsilon(x)\right)^{n} \\
  &   =  & C \supmasyL^{n+1+d/\alpha}(n+1)^{-d/\alpha}
           \left[\frac{1}{n+1}\left(1-\left(1-\frac{b_{\varepsilon}(x)}{\supmasyL}\right)^n\right)\right. \\
  &      & \left. - \frac{b_\varepsilon(y)}{\supmasyL}\left(1+\frac{1}{n}\right)^{d/\alpha}\frac{a\eta^{-\alpha}}{A(n+1)}
           +\left(1+\frac{1}{n}\right)^{d/\alpha}\right. \\
  &      & \left.+\left(1-\frac{b_\varepsilon(x)}{\supmasyL}\right)^n\left(n+1\right)^{d/\alpha}\right] \\
  & \leq & C \supmasyL^{n+1+d/\alpha}(n+1)^{-d/\alpha}\left[\frac{2}{n+1}+\left(1+\frac{1}{n}\right)^{d/\alpha}
          \left(1-\frac{\eta^{-\alpha}(a/A)^2}{n+1}\right)\right] \\
  & \leq & C \supmasyL^{n+1+d/\alpha}(n+1)^{-d/\alpha}\left[\frac{2}{n+1} +
           \left(1+\frac{p}{n}\right)
          \left(1-\frac{\eta^{-\alpha}(a/A)^2}{n+1}\right)\right] \\
  & \leq &  C \supmasyL^{n+1+d/\alpha}(n+1)^{-d/\alpha}\left[1 - 
            \frac{1}{n+1}\left( \eta^{-\alpha}(a/A)^2-2 - 2p\right)
          \right] \\
  &   =  & C \supmasyL^{n+1+d/\alpha}(n+1)^{-d/\alpha},
\end{eqnarray*}
where the second inequality follows from (\ref{eq:supbx}), (A.4) and (\ref{eq:hilfe}).
\qed

Using the above lemmas we may estimate $\Gamma_\varepsilon^n$ and in consequence also the exponent 
operator $e^{t\Aapprox{\varepsilon}}=e^{-t\supmasyL}e^{t\Gamma_\varepsilon}$.

\begin{lemma}\label{lm:SE1} Assume (A.1) -- (A.4). There exists a constant $C$ such that for all $x\in\R$ and all nonnegative 
$\varphi\in\Bbound$ such that $x\notin \supp(\varphi)$ we have
  \begin{equation}
 e^{t\Aapprox{\varepsilon}}\varphi(x) \leq  C t \int \varphi(z) |z-x|^{-\alpha-d}\,dz,\quad \varepsilon>0.
\end{equation}
\end{lemma}
\dowod
By Lemma \ref{lm:GammaDensity} and Lemma \ref{lm:Estimate1} for every $\varphi$ such that $x\not\in\supp(\varphi)$ we get
$$
 \Gamma_\varepsilon^n\varphi(x) \leq \int \varphi(y) C n \supmasyL^{n-1}|y-x|^{-\alpha-d}\,dy,
$$
and
\begin{eqnarray*}
 e^{t\Aapprox{\varepsilon}}\varphi(x) 
 & \leq & Ce^{-t\supmasyL} \sum_{n=1}^\infty 
          \frac{t^n n \supmasyL^{n-1}}{n!}  \int \varphi(y) |y-x|^{-\alpha-d}\,dy \\
 &   =  & Ce^{-t\supmasyL} t \sum_{n=0}^\infty 
          \frac{t^n\supmasyL^{n}}{n!}  \int \varphi(y) |y-x|^{-\alpha-d}\,dy \\
 &   =  & C t \int \varphi(y) |y-x|^{-\alpha-d}\,dy.
\end{eqnarray*}
\qed

The following lemma seems to be known. We include here the proof for reader's convenience.

\begin{lemma} For every $p\in [0,\infty)$ exists a constant $C=C(p)$ such that
\begin{equation}\label{eq:Szereg}
  \sum_{n=1}^\infty \frac{x^{n+p}}{n!n^p}\leq C (e^{x}-1),\quad x>0.
  \end{equation}
\end{lemma}
\dowod For $p\in [0,1]$ we have by Jensen's inequality
\begin{eqnarray*}
  \sum_{n=1}^\infty \frac{(x/n)^px^n}{n!} 
  & \leq & \frac{e^x-1}{(e^x-1)^p}\left(\sum_{n=1}^\infty \frac{x^{n+1}}{n\cdot n!} \right)^p \\
  &  =   & \frac{e^x-1}{(e^x-1)^p}\left(\sum_{n=1}^\infty \frac{x^{n+1}}{(n+1)!} \frac{n+1}{n}\right)^p \\
  & \leq &  2^p (e^x-1),\quad x>0.
\end{eqnarray*}
For $p\geq 1$ we have
$$
  \sum_{n=1}^\infty \frac{x^{n+p}}{n!n^p} 
  =  \sum_{n=1}^\infty \int_0^x \frac{u^{n+p-1}}{n!n^{p-1}}\,du\frac{n+p}{n} \leq (p+1) \int_0^x \sum_{n=1}^\infty \frac{u^{n+p-1}}{n!n^{p-1}}\,du,
$$
and the Lemma follows by induction.
\qed

\begin{lemma}\label{lm:SE2} Assume (A.1), (A.3) and (A.4). Then there exists $C$ such that for every nonnegative $\varphi\in \Bbound\cap L_1(\R)$ we have
  \begin{equation}
 e^{t\Aapprox{\varepsilon}}\varphi(x) \leq C t^{-d/\alpha} \int \varphi(y)\,dy + e^{-tb_\varepsilon(x)}\varphi(x),
 \quad x\in\R,\,\varepsilon>0,t>0. 
\end{equation}
\end{lemma}
\dowod By Lemma \ref{lm:Estimate3} for every $\varphi\in \Bbound\cap L_1(\R)$ we get
$$
 \Gamma_\varepsilon^n\varphi(x) \leq  c (\supmasyL)^{n+d/\alpha}
    n^{-d/\alpha}  \int \varphi(y)\,dy + \left(\supmasyL-b_\varepsilon(x)\right)^{n}\varphi(x),
$$
and by (\ref{eq:Szereg}) we obtain
\begin{eqnarray*}
 e^{t\Aapprox{\varepsilon}}\varphi(x) 
 & \leq & e^{-t\supmasyL} \left[c \int \varphi(y)\,dy \sum_{n=1}^\infty 
        \frac{ t^n \supmasyL^{n+d/\alpha}}{n!n^{d/\alpha}} + 
        e^{t\left(\supmasyL-b_\varepsilon(x)\right)}\varphi(x)\right]         \\
 & \leq & c t^{-d/\alpha} \int \varphi(y)\,dy +  e^{-tb_\varepsilon(x)}\varphi(x). 
 \end{eqnarray*}
 \qed

\dowod of Theorem \ref{th:MainT}. Let $t>0$, $\varphi\in\Bbound$, and $x\in\R$. 
Using Lemma \ref{lm:SE1} for $\indyk{B(x,t^{1/\alpha})^c}\varphi$ and Lemma \ref{lm:SE2} for $\indyk{B(x,t^{1/\alpha})}\varphi$ we obtain
\begin{eqnarray*}
  e^{t\Aapprox{\varepsilon}}\varphi(x)
  &  =   & e^{t\Aapprox{\varepsilon}}[\indyk{B(x,t^{1/\alpha})^c}\varphi](x) + e^{t\Aapprox{\varepsilon}}[\indyk{B(x,t^{1/\alpha})}\varphi](x) \\
  & \leq & c\left[\int_{B(x,t^{1/\alpha})^c} \varphi(y) t|y-x|^{-\alpha-d}\, dy +
            \int_{B(x,t^{1/\alpha})} \varphi(y) t^{-d/\alpha}\,dy \right] \\
  &      & + \,e^{-tb_\varepsilon(x)}\varphi(x) \\
  & \leq & c\int \varphi(y) \min(t^{-d/\alpha},t|y-x|^{-\alpha-d})\, dy + e^{-tb_\varepsilon(x)}\varphi(x)
\end{eqnarray*}
\qed

We prove now that the continuity of the jump intensity $f(x,y)$ yields a regularity of $\Aapprox{\varepsilon}$.

\begin{lemma}\label{l:Aepsilon} 
  If (A.5) holds then for all 
  $\varepsilon>0$ the operator $\Aapprox{\varepsilon}$ maps $\Cvanish$ into $\Cvanish$.
\end{lemma}
\dowod Let $\varphi\in\Cvanish$ and $\eta>0$. We choose $r>0$ such that $|\varphi(y)|\leq \frac{\eta\varepsilon^\alpha}{4A}$ for $|y|\geq r$, and 
using (A.1) and (\ref{eq:supbx}), for $|x|>r$ we get
\begin{eqnarray*}
  |\Aapprox{\varepsilon}\varphi(x)| 
  & \leq & \int_{|y|<r} |\varphi(y)|\, f_{\varepsilon}(x,y)\,dy + \int_{|y|\geq r} |\varphi(y)|\, f_{\varepsilon}(x,y)\,dy
  + |\varphi(x)| \int_{\R}  f_{\varepsilon}(x,y)\,dy \\
  & \leq & M \|\varphi\|_\infty \int_{|y|<r} |y-x|^{-\alpha-d}\, dy + \frac{\eta\varepsilon^\alpha}{4A} 2 b_{\varepsilon}(x) \\
  & \leq & c_1 \|\varphi\|_\infty  (|x|-r)^{-\alpha-d} r^d + \frac{\eta\varepsilon^\alpha}{4A} 2A\varepsilon^{-\alpha}.
\end{eqnarray*}
For $|x|>r+(2c_1r^d\|\varphi\|_\infty\eta^{-1})^{1/(\alpha+d)}$ we obtain $|\Aapprox{\varepsilon}\varphi(x)|< \eta$. This yields 
$ \lim_{|x|\to\infty} \Aapprox{\varepsilon}\varphi(x)=0$.

Let $x_0\in\R$. For $|x-x_0|<\varepsilon/2$, and all $y\in\R$ we have 
\begin{eqnarray*}
  f_{\varepsilon}(x,y)
  & \leq & M \left[\max(\varepsilon,|y-x|)\right]^{-\alpha-d}
  \leq  M \left[\max(\varepsilon,|y-x_0|-\varepsilon/2)\right]^{-\alpha-d} \\
  & \leq & M2^{\alpha+d}\left(|y-x_0|+\varepsilon/2\right)^{-\alpha-d}.
\end{eqnarray*}
By the continuity of $x\to f(x,y)$ and the dominated convergence theorem we get
$$
  \lim_{x\to x_0}\int |f_{\varepsilon}(x,y)-f_{\varepsilon}(x_0,y)|\, dy =0.
$$
It follows that $b_\varepsilon(x)$ is continuous on $\R$ for each $\varepsilon>0$. Consequently for every
$\varphi\in \Cvanish$ 
we have
\begin{eqnarray*}
  |\Aapprox{\varepsilon}\varphi(x)-\Aapprox{\varepsilon}\varphi(x_0)| 
  & \leq & \int |\varphi(y)|\, |f_{\varepsilon}(x,y)-f_{\varepsilon}(x_0,y)|\,dy \\
  &      &  + |\varphi(x)b_\varepsilon(x)-\varphi(x_0)b_\varepsilon(x_0)|\,\to 0,\,\mbox{ as } x\to x_0.
\end{eqnarray*}
Thus $\Aapprox{\varepsilon}\varphi\in\Cvanish$.
\qed

We note that the operators $\Aapprox{\varepsilon}$ approximate $\gener$.

\begin{lemma}\label{l:AaA} If (A.2) is satisfied then for all $\varphi\in C^2_c(\R)$ we have
$$\lim_{\varepsilon\to 0}\|\gener\varphi - \Aapprox{\varepsilon}\varphi\|_\infty=0$$
\end{lemma}

\dowod For all $\varphi\in C^2_c(\R)$, $x\in\R$, and $\varepsilon\in(0,1)$ we have
\begin{eqnarray*}
\Aapprox{\varepsilon}\varphi(x) & =  &
\int \left(\varphi(y)-\varphi(x)\right)f_{\varepsilon}(x,y)\,dy\\
&  =   & \int_{\R} \left(\varphi(x+h)-\varphi(x)-\indyk{B(0,1)}(h)\scalp{h}{\nabla\varphi(x)}\right)
          f_{\varepsilon}(x,x+h)\,dh \\
&      & + \frac{1}{2}\int_{|h|<1} \scalp{h}{\nabla\varphi(x)} \left(f_{\varepsilon}(x,x+h)-f_{\varepsilon}(x,x-h) \right)\,dh.
\end{eqnarray*}
Furthermore, it follows from (A.1) and (A.2) that
\begin{equation}\label{eq:intf}
  \int_{|h|<\varepsilon} |h||f(x,x+h)-f(x,x-h)|\,dh \leq c  \varepsilon^{1-\alpha}\max(1-\alpha, 0).
\end{equation}
Using Taylor's expansion for $\varphi$ we obtain
\begin{eqnarray*}
  |\gener\varphi(x) - \Aapprox{\varepsilon}\varphi(x)| 
  & \leq &  \int_{|h|<\varepsilon} \left|\varphi(x+h)-\varphi(x)-\scalp{h}{\nabla\varphi(x)}\right|
            f(x,x+h)\,dh \\
  &      & +\, \frac{1}{2}|\nabla\varphi(x)|\int_{|h|<\varepsilon} |h| \left|f(x,x+h)-f(x,x-h) \right|\,dh \\
  & \leq & c\sup_{y\in\R,|\beta|\leq 2}  |\partial^\beta\varphi(y)|\left(\int_{|h|<\varepsilon} |h|^{2-\alpha-d}\,dh\right. \\
  &      & \left. +\, \int_{|h|<\varepsilon} |h| \left|f(x,x+h)-f(x,x-h) \right|\,dh\right)\\
  & \leq & c\sup_{y\in\R,|\beta|\leq 2}  |\partial^\beta\varphi(y)|\left(\varepsilon^{2-\alpha}+ \varepsilon^{1-\alpha}\max(1-\alpha, 0)\right),
\end{eqnarray*}
which yields $\lim_{\varepsilon\to 0}\|\gener\varphi - \Aapprox{\varepsilon}\varphi\|_\infty=0$.
\qed

We are ready now to prove Lemma \ref{l:AonCvanish}.

\dowod {\it of Lemma \ref{l:AonCvanish}}. The continuity of $\gener\varphi$ for $\varphi\in C^2_c(\R)$ follows from
Lemma \ref{l:Aepsilon} and Lemma \ref{l:AaA}. For $\varphi\in C^2_c(\R)$ and $\eta>0$ we choose $\varepsilon>0$ such that
$\|\gener\varphi - \Aapprox{\varepsilon}\varphi\|_\infty<\eta/2$ and $r>0$ such that $|\Aapprox{\varepsilon}\varphi(x)|<\eta/2$ for
all $|x|>r$.
We obtain
$
  |\gener\varphi(x)| \leq \|\gener\varphi - \Aapprox{\varepsilon}\varphi\|_\infty + |\Aapprox{\varepsilon}\varphi(x)| < \eta,
$
for all $|x|>r$.
\qed

\dowod {\it of Theorem \ref{th:main}}. By Lemma \ref{l:AaA} we have
$$ 
 \lim_{\varepsilon\to 0}\|\gener\varphi - \Aapprox{\varepsilon}\varphi\|_\infty=0
$$
for every $\varphi\in C_\infty^2(\R)$. 
A closure of $\gener$ is a generator of a semigroup and from the Hille-Yosida theorem it follows that the range of $\lambda-\gener$ is dense in $\Cvanish$ 
and therefore by Theorem 5.2 in \cite{Trotter58} (see also \cite{Hasegawa}) we get
$$
  \lim_{\varepsilon\downarrow 0} \|e^{t\Aapprox{\varepsilon}}\varphi-P_t\varphi\|_\infty =0,
$$
for every $\varphi\in \Cvanish$. By Theorem \ref{th:MainT} this yields
$$
  P_t\varphi(x) \leq c_1  \int 
  \varphi(z) \min\left(t^{-d/\alpha},\frac{t}{|z-x|^{\alpha+d}}\right) dz,
$$
for every nonnegative $\varphi\in \Cvanish$.
\qed

\end{document}